\documentclass{elsarticle}

\usepackage{amssymb}
\usepackage{amsthm}
\usepackage{amsmath}
\usepackage{enumitem}
\usepackage{nicefrac}
\usepackage{hyperref}

\DeclareMathOperator{\RCAn}{RCA_0}
\DeclareMathOperator{\WKLn}{WKL_0}
\DeclareMathOperator{\WKL}{WKL}
\DeclareMathOperator{\ACAn}{ACA_0}
\DeclareMathOperator{\ACA}{ACA}
\DeclareMathOperator{\ATRn}{ATR_0}
\DeclareMathOperator{\WETR}{WETR}
\DeclareMathOperator{\SETR}{SETR}
\DeclareMathOperator{\WO}{WO}
\newcommand{\N}{\mathbb{N}}
\newcommand{\period}{\text{.}}
\newcommand{\komma}{\text{,}}
\newcommand{\restr}[2]{#1\!\upharpoonright\!#2}

\numberwithin{equation}{section}

\begin{document}
	
\theoremstyle{plain}
\newtheorem{lemma}{Lemma}
\newtheorem{proposition}[lemma]{Proposition}
\newtheorem{corollary}[lemma]{Corollary}
\newtheorem{theorem}[lemma]{Theorem}
	
\theoremstyle{definition}
\newtheorem{definition}[lemma]{Definition}
\newtheorem{example}[lemma]{Example}

\theoremstyle{remark}
\newtheorem{remark}[lemma]{Remark}

\begin{frontmatter}
	\title{Weak and Strong Versions of\\Effective Transfinite Recursion\tnoteref{t1}}
	\tnotetext[t1]{\copyright{} 2022. This manuscript version is made available under the CC-BY-NC-ND 4.0 license \url{https://creativecommons.org/licenses/by-nc-nd/4.0}. It is the accepted version of a paper	published in the Annals of Pure and Applied Logic 174(4) 2023, article no.~103232, 15~pp. Funded by the Deutsche Forschungsgemeinschaft (DFG, German Research Foundation) -- Project number 460597863.}
	\author[P. Uftring]{Patrick Uftring}
	\address{Department of Mathematics, Technical University of Darmstadt, Schloss\-garten\-str.~7, 64289~Darmstadt, Germany}
	\ead{uftring@mathematik.tu-darmstadt.de}
	\begin{abstract}
		Working in the context of reverse mathematics, we give a fine-grained characterization result on the strength of two possible definitions for Effective Transfinite Recursion used in literature. Moreover, we show that $\Pi^0_2$-induction along a well-order $X$ is equivalent to the statement that the exponentiation of any well-order to the power of $X$ is well-founded.
	\end{abstract}
	\date{11. October, 2022}
\end{frontmatter}

\section{Introduction}

The principle of \emph{Arithmetical Transfinite Recursion} allows us to construct a set~$Y$ along well-orders $X$, where each stage ${Y_x := \{n \in \N \text{ with } (x, n) \in Y\}}$ is defined using an arithmetical formula $\varphi(n, x, Z)$ such that a pair $(x, n)$ enters $Y$ if and only if $\varphi(n, x, Y^x)$ holds. Here, $Y^x$ denotes all pairs of previous stages, i.e., the set $\{(y, m) \in Y \mid y <_X x \text{ and } m \in \N\}$. More on this topic can, e.g., be found in Section V.2 of \cite{Simpson}.

In this paper, we study different definitions of \emph{Effective Transfinite Recursion} in the context of reverse mathematics. Inspired by the work due to Church, Kleene, and Rogers (cf.~\cite[Section I.3]{Sacks90}), \emph{Weak} Effective Transfinite Recursion ($\WETR$) is defined like Arithmetical Transfinite Recursion but where $\varphi$ is required to be a $\Delta^0_1$-formula (cf.~\cite[Definition 6.12]{Dzhafarov}).
In \cite{GreenbergMontalban}, Greenberg and Montalb\'{a}n applied $\WETR$ in order to construct several embeddings between linear orders (cf.~\cite[Lemma~2.4 and Lemma~3.27]{GreenbergMontalban}).
In \cite{Dzhafarov}, Dzhafarov et al.~showed that $\WETR$ can be proved in $\ACAn$. They used this result to show that a certain statement about Borel codes and evaluation maps already implies the much stronger system of $\ATRn$ (cf.~\cite[Proposition 6.13 and Theorem 6.14]{Dzhafarov}).

In Section \ref{sec:WETR} of this paper, we will not only show that $\WETR$ is equivalent to Weak K\H{o}nig's Lemma, but also that each instance of $\WETR_X$ for a fixed well-order~$X$ is equivalent to the disjunction of $\WKL$ and $\Pi^0_2$-induction along~$X$. Moreover, we present a rule version of $\WETR$ that is already admissible over $\RCAn$.

In \cite{Freund20}, Freund introduced a stronger variant of this principle of Weak Effective Transfinite Recursion. \emph{Strong} Effective Transfinite Recursion ($\SETR$) is defined like $\WETR$ but $\varphi(n, x, Z)$ is $\Sigma^0_1$ and only needs to be $\Delta^0_1$ (witnessed by some fixed $\Pi^0_1$-formula $\lnot \psi(n, x, Z)$) for $Z$ that correspond to a set that is recursively defined using $\varphi$ along $\restr{X}{x} := \{y \in X \mid y <_X x\}$. Adapting the proof by Dzhafarov et al., Freund showed that this variant is also implied by $\ACAn$. Later, in \cite{Freund21}, he used $\SETR$ in order to prove a type-$2$ analogue of the statement that every normal function has a fixed point.

In Section \ref{sec:SETR}, we will show that this new principle $\SETR$ is not only implied by but also equivalent to $\ACAn$. Moreover, we will prove that each instance $\SETR_X$ for a fixed well-order $X$ is equivalent to $\Pi^0_2$-induction along $X$ and to the statement that exponentiation to the power of $X$, a certain map between linear orders, preserves well-orders. Together with the results from before, we will see that $\WETR_X$ is equivalent to the disjunction of $\WKL$ and $\SETR_X$.

\section{Strong Effective Transfinite Recursion}\label{sec:SETR}

The following definitions are taken from \cite{Freund20}. We begin this section by recalling the definition of a \emph{recursively defined family}:

\begin{definition}[Recursively defined family]\label{def:recursively_defined_faimly}
	Given a well-order $X$ together with a formula $\varphi(n, x, Z)$, we express that $Y$ is the \emph{recursively defined family} with respect to $X$ and $\varphi$ using the following formula $H_\varphi(X, Y)$:
	\begin{equation*}
		H_{\varphi}(X, Y) :\Longleftrightarrow Y = \{(x, n) \mid x \in X \text{ and } \varphi(n, x, Y^x)\}\komma
	\end{equation*}
	where $Y^x$ is the restriction of $Y$ to pairs $(x', n)$ with $x' <_X x$.
\end{definition}

We then define \emph{Strong Effective Transfinite Recursion} similar to Arithmetical Transfinite Recursion but where we only consider formulas $\varphi(n, x, Z)$ that are decidable in the case where $Z$ corresponds to earlier stages of the principle:

\begin{definition}[Strong Effective Transfinite Recursion]\label{def:SETR}
	Given a well-order $X$, we define the principle $\SETR_X$ of \emph{Strong Effective Transfinite Recursion} along $X$:
	Consider $\Sigma^0_1$-formulas $\varphi(n, x, Z)$ and $\psi(n, x, Z)$. If we have
	\begin{equation}\label{eq:assumption_SETR}
		\forall n \in \N\ \forall x \in X\ \forall Z \subseteq \N \ \big(H_{\varphi}(\restr{X}{x}, Z) \to \big(\varphi(n, x, Z) \leftrightarrow \lnot \psi(n, x, Z)\big)\big)\komma
	\end{equation}
	where $\restr{X}{x}$ denotes the suborder of $X$ that only includes elements $x' \in X$ with $x' <_X x$, then there exists some set $Y$ with $H_{\varphi}(X, Y)$.
	We write $\SETR$ for the principle that $\SETR_X$ holds for any well-order~$X$.
\end{definition}

\begin{remark}\label{rem:single_formula}
	We can express the principle $\SETR_X$ using a single formula. For this, we use Kleene's normal form (cf.~\cite[Theorem II.2.7]{Simpson}) and $\Delta^0_1$-comprehension in order to write any two $\Sigma^0_1$-formulas $\varphi(n, x, Z)$ and $\psi(n, x, Z)$ equivalently as
	\begin{equation*}
		\varphi'(n, x, Z) := \exists m\ (n, x, Z[m]) \in U \text{ and } \psi'(n, x, Z) := \exists m\ (n, x, Z[m]) \in V\komma
	\end{equation*}
	respectively, for sets $U$ and $V$. Here and later in this paper, $Z[m]$ denotes the initial segment of $Z$ of length $m$, i.e., the unique finite $0$-$1$-sequence of length $m$ with $Z[m]_i = 1 \leftrightarrow i \in Z$ for all $i <_\N m$. The whole principle $\SETR_X$, therefore, already holds if we only assume $\SETR_X$ for the formulas $\varphi'(n, x, Z)$ and $\psi'(n, x, Z)$.
	Finally, $\SETR$ is then expressed as the formula 
	\begin{equation*}
		\SETR := \forall X \big(\WO(X) \to \SETR_X\big)\komma
	\end{equation*}
	where $\WO(X)$ formulates that $X$ is a well-order.
\end{remark}

Next, we recall the definition for exponentiation of well-orders from \cite[Definition~2.1]{Hirst94} where we additionally allow all linear orders as base:

\begin{definition}
	Given a linear order $\alpha$ together with a least element $0_\alpha \in \alpha$ (if it exists), and a well-order $\beta$, we define the set $\alpha^\beta$ to contain all elements of the form
	\begin{equation*}
		\langle (b_0, a_0), \dots, (b_{n-1}, a_{n-1}) \rangle
	\end{equation*}
	for $n \in \N$, $a_0, \dots, a_{n-1} \in \alpha$ and $b_0, \dots, b_{n-1} \in \beta$ with $a_i \neq 0_\alpha$ for all $i <_\N n$ and $b_i >_\beta b_{i+1}$ for all $i <_\N n - 1$. The set $\alpha^\beta$ is ordered lexicographically, i.e., given two elements $\sigma, \tau \in \alpha^\beta$, we have $\sigma <_{\alpha^\beta} \tau$ if and only if one of the following holds:
	\begin{itemize}
		\item $\tau$ strictly extends $\sigma$, or
		\item for the least index $i$ with $\sigma_i \neq \tau_i$, we have either $b <_\beta b'$ or $b = b'$ and $a <_\alpha a'$ for $(b, a) := \sigma_i$ and $(b', a') := \tau_i$.
	\end{itemize}
\end{definition}

In contrast to our source, we have less strict requirements on $\alpha$ and $\beta$. We do this in order to be able to map any linear order $\alpha$ to $\alpha^\beta$ for any well-order $\beta$. Moreover, during the proof of Lemma \ref{lem:varphi'_terminates}, we will even need to allow partial-orders~$\alpha$. This is not a problem since the definition neither requires $\alpha$ to be well-founded nor to be a linear order.
Now, we have collected all necessary definitions in order to state the main result of this section:

\begin{theorem}[$\RCAn$]\label{thm:SETR}
	Consider some well-order $X$. Then, the following are equivalent:
	\begin{enumerate}[label=\alph*)]
		\item Strong Effective Transfinite Recursion along $X$,
		\item Transfinite $\Pi^0_2$-induction along $X$,
		\item The map $\alpha \mapsto \alpha^X$ preserves well-orders.
	\end{enumerate}
\end{theorem}

We divide the directions for the proof of this theorem into three propositions. We begin by using $\SETR_X$ to construct a set $Y$ such that, under the assumption $H_\varphi(X, Y)$, we have $\Pi^0_2$-induction along $X$.

\begin{proposition}[$\RCAn$]\label{prop:SETR_X_implies_induction}
	Consider some well-order $X$ and assume $\SETR_X$. Then, we can show transfinite $\Pi^0_2$-induction along $X$.
\end{proposition}

\begin{proof}
	Let $P(x, n, m)$ be an arbitrary $\Delta^0_0$-formula such that the $\Pi^0_2$-formula
	\begin{equation*}
		Q(x) := \forall n\ \exists m\ P(x, n, m)
	\end{equation*}
	is progressive, i.e., we have $\forall y <_X x\ Q(y) \to Q(x)$ for any $x \in X$. Using $\SETR_X$, we want to show that $Q(x)$ holds for all $x \in X$.
	
	For this, we define the following $\Sigma^0_1$-formulas $\varphi(n, x, Z)$ and $\psi(n, x, Z)$:
	\begin{align*}
		\varphi(n, x, Z) &:= \exists m\ \big(P(x, n, m) \land \forall (y, n') \leq_\N m\ \big(y <_X x \to (y, n') \in Z\big)\big)\\
		\lnot\psi(n, x, Z) &:= \forall y <_X x\ \forall n'\ \big(\forall m <_\N (y, n')\ \lnot P(x, n, m) \to (y, n') \in Z\big)\period
	\end{align*}
	For every $x \in X$ and $Z \subseteq \N$, we want to show that under the assumption ${H_\varphi(\restr{X}{x}, Z)}$ the equivalence $\varphi(n, x, Z) \leftrightarrow \lnot \psi(n, x, Z)$ holds for all $n \in \N$.
	
	``$\rightarrow$'': Assume that $\varphi(n, x, Z)$ holds and is witnessed by some $m \in \N$. Consider some $y <_X x$ and $n' \in \N$ such that $\lnot P(x, n, m')$ holds for all $m' <_\N (y, n')$. Each such $m'$ must be different from $m$. We conclude $(y, n') \leq_\N m$. With $y <_X x$ and the assumption on $m$, we conclude $(y, n') \in Z$.
	
	``$\leftarrow$'': Assume that $\varphi(n, x, Z)$ does not hold. We may assume that there is some pair $(y, n') \notin Z$ with $y <_X x$. Otherwise, using the assumption ${H_\varphi(\restr{X}{x}, Z)}$, we have $\varphi(n', y, Z)$ for all $n' \in \N$ and $y <_X x$. This implies $Q(y)$ for all $y <_X x$ and, therefore, $Q(x)$ since $Q$ is progressive. Finally, together with the assumption $(y, n') \in Z$ for all $y <_X x$ and $n' \in \N$ this produces the contradiction $\varphi(n, x, Z)$. Now that there is some pair $(y, n') \notin Z$ with $y <_X x$, we may assume that it is the pair with the smallest code that experiences this property. Since $\varphi(n, x, Z)$ does not hold, we conclude that any $m \in \N$ with $P(x, n, m)$ must be greater or equal to the code of our smallest $(y, n')$. Consequently, we have $\lnot P(x, n, m')$ for all $m' <_\N (y, n')$. Finally, together with $(y, n') \notin Z$ this means that $\lnot \psi(n, x, Z)$ cannot hold.
	
	We can now apply $\SETR_X$, which produces some set $Y$ with $H_\varphi(X, Y)$. Assume that there is some pair $(x, n) \notin Y$ for $x \in X$ and $n \in \N$. Let it be the smallest element in $(X \times \N) \setminus Y$ having this property with respect to the lexicographical order on $X \times \N$. This is a well-order since $X$ is one. With $(x, n)$ being the smallest such element, $(y, n') \in Y$ holds for all $y <_X x$ and $n' \in \N$. Moreover, with $H_\varphi(X, Y)$, this means that $Q(y)$ holds for all $y <_X x$. Therefore, we have $Q(x)$ since $Q$ is progressive. By definition of $\varphi$ and the fact that all pairs $(y, n)$ with $y <_X x$ satisfy $(y, n) \in Y$, this implies the contradiction $\varphi(n, x, Y^x)$. We conclude that all pairs $(x, n) \in X \times \N$ are elements of $Y$. Finally, by $H_\varphi(X, Y)$ and the definition of $\varphi$, this entails that $Q(x)$ must hold for all~${x \in X}$.
\end{proof}

With this, we can already prove the equivalence of the first two statements in our theorem, although this will also follow as soon as we have proven the remaining two propositions:

\begin{corollary}[$\RCAn$]\label{cor:induction_to_SETR_X}
	For any well-order $X$, the principles of $\SETR_X$ and $\Pi^0_2$-induction along $X$ are equivalent. Moreover, $\SETR$ is equivalent to $\ACA_0$.
\end{corollary}

\begin{proof}
	The direction from $\SETR_X$ to $\Pi^0_2$-induction along $X$ is covered by the previous proposition. The other direction will follow by combining Propositions \ref{prop:induction_proves_alpha_x}~and~\ref{prop:alpha_x_proves_SETR_X}. Another way of proving this direction, that relies on established literature, is the following:
	Inspecting Freund's proof that $ACA_0$ implies $\SETR$ (cf.~\cite[Theorem~5]{Freund20}), that he adapted from the similar result that $ACA_0$ implies $\WETR$ due to Dzhafarov et al. (cf.~\cite[Proposition~6.13]{Dzhafarov}), we see that already $\Pi^0_2$-induction along $X$ implies $\SETR_X$. Note for this that Freund applies transfinite induction along $X$ to the conjunction of two statements $(I)$ and $(II)$. This conjunction is a $\Pi^0_2$-formula.
	
	For the equivalence of $\SETR$ and $\ACA_0$, we apply the result that transfinite induction for $\Sigma^0_1$-formulas and for arithmetical formulas is equivalent to arithmetical comprehension (cf.~\cite[Corollary~3]{Hirst99}).
\end{proof}

Using Proposition \ref{prop:SETR_X_implies_induction}, we can also begin to transport results about $\Pi^0_2$-induction to Strong Effective Transfinite Recursion:

\begin{corollary}
	$\RCAn$ cannot prove $\SETR_\omega$ but it can prove the equivalence of this principle to $\forall n\ \SETR_n$.
\end{corollary}

\begin{proof}
	It is clear that $\SETR_\omega$ implies $\SETR_n$ for all $n \in \N$. Assume that $\forall n\ \SETR_n$ holds. Using Proposition \ref{prop:SETR_X_implies_induction}, we have $\Pi^0_2$-induction up to any natural number $n$. Therefore, we have $\Pi^0_2$-induction along $\N$. Finally, by Corollary \ref{cor:induction_to_SETR_X}, this yields $\SETR_\omega$.
	$\Pi^0_2$-induction, however, is not provable in $\RCAn$ (cf.~\cite{Parsons} also \cite[Propositions~4 and 7]{KP78}).
\end{proof}

We continue with the second direction for the proof of Theorem \ref{thm:SETR}:

\begin{proposition}[$\RCAn$]\label{prop:induction_proves_alpha_x}
	Consider some well-order $X$. $\Pi^0_2$-induction along $X$ implies that $\alpha \mapsto \alpha^X$ preserves well-orders.
\end{proposition}

\begin{proof}
	For contradiction, assume that we can find some infinite descending sequence~$(f_n)_{n \in \N}$ in $\alpha^X$ for some well-order $\alpha$. Given an element $g \in \alpha^X$, we write $g(x)$ for the unique value $a \in \alpha$ with $(x, a) \in g$ if it exists or $0_\alpha$ if it does not.
	For any $x \in X$, we have the sequence $(f^x_n)_{n \in \N}$ where $f^x_n \in \alpha^X$ is defined like the finite sequence $f_n$ but where any pair $(y, a) \in f_n$ with $y <_X x$ is removed, for all~$n \in \N$. For example, if $0_X$ is the smallest element of $X$, then we have $f^{0_X}_n = f_n$ for all~$n \in \N$.
	We define the following $\Pi^0_2$-formula:
	\begin{equation*}
		\varphi(x) := \text{ For all $m \in \N$, there is some $n \in \N$ with $f^x_n >_\N m$}
	\end{equation*}
	Given an element $x \in X$, this formula $\varphi(x)$ holds if and only if~$(f^x_n)_{n \in \N}$ reaches infinitely many distinct values. Let $x \in X$ be the greatest element of $X$ satisfying~${f_0(x) \neq 0_\alpha}$. By definition of $\alpha^X$ and the assumption that $(f_n)_{n \in \N}$ is descending, we see that for any $n \in \N$ and $y \in X$ with $f_n(y) \neq 0_\alpha$, we have $y \leq_X x$. Thus, $\varphi(x)$ cannot hold since otherwise $(f_n(x))_{n \in \N}$ would be a sequence in the well-order~$\alpha$ that never ascends but descends infinitely often. Using $\Pi^0_2$-induction along $X$, we conclude that there must be some least $x \in X$ such that $\varphi(x)$ does not hold. In the following case distinction, we show that this leads to a contradiction:
	
	\begin{itemize}
		\item Base:
		
		Let $x := 0_X$, i.e., $x$ is the least element of $X$. As already mentioned in the introduction of this proof, we have $f^x_n = f_n$ for all $n \in \N$. Thus, $(f^x_n)_{n \in \N}$ reaches infinitely many different values, which implies $\varphi(x)$.
		\item Successor:
		
		Assume that there exists some $y \in X$ such that $x$ is the least element in~$X$ with $x >_X y$. Since $(f_n)_{n \in \N}$ is decreasing, this implies that $(f^x_n)_{n \in \N}$ is weakly decreasing. Together with the fact that this sequence only reaches finitely many distinct values, we know that it must be constant starting from some index $m \in \N$. Therefore, since $(f^y_n)_{n \in \N}$ consists of infinitely many values, we know that $(f_n(y))_{n \in \N}$ must have the same property. Also, since $(f_n)_{n \in \N}$ is decreasing, we know that $(f^y_n)_{n \in \N}$ is weakly decreasing. Combining this with the facts that $(f^x_n)_{n \geq_\N m}$ is constant and that $x$ is the successor of $y$, we derive that $(f_n(y))_{n \geq_\N m}$ is weakly decreasing in $\alpha$. In summary, $(f_n(y))_{n \geq_\N m}$ is a weakly decreasing sequence in $\alpha$ consisting of infinitely many values. This is impossible in a well-order.
		\item Limit:
		
		Assume that $x$ is neither a successor nor equal to $0_X$. Using the same argument as before, we know of some $m \in \N$ such that $(f^x_n)_{n \geq_\N m}$ is constant. Let $y \in X$ be some value in $X$ with $x >_X y$ and $y >_X z$ for all $z \in X$ with $f_m(z) \neq 0_\alpha$. This value $y$ exists since $x$ is neither a successor nor equal to $0_X$. We have $f^x_m = f^y_m$. By induction hypothesis, $(f^y_n)_{n \geq_\N m}$ reaches infinitely many different values. Therefore, there must be some index $n >_\N m$ with $f^y_n \neq f^y_m$. Since $f^x_n = f^x_m$, we conclude that there are distinct finite sequences $s$ and $t$ with $f^y_n = f^x_n * s$ and $f^y_m = f^x_m * t$ where the operator ``$*$'' concatenates sequences. From $f^x_m = f^y_m$, we know that $t$ is empty. Since $s$ is distinct from $t$, we have $f^y_n >_{\alpha^X} f^y_m$. This implies $f_n >_{\alpha^X} f_m$ by definition of the order on~$\alpha^X$. Finally, this inequality and $n >_\N m$ contradict our assumption that $(f_n)_{n \in \N}$ was decreasing.
	\end{itemize}
\end{proof}

We finish with the third and last direction of Theorem \ref{thm:SETR}:

\begin{proposition}[$\RCAn$]\label{prop:alpha_x_proves_SETR_X}
	Consider some well-order $X$. If $\alpha \mapsto \alpha^X$ preserves well-orders, then $\SETR_X$ holds.
\end{proposition}

The proof for this proposition is more involved than those of the previous two. Because of this, we split it into Definition \ref{def:varphi'}, and Lemmas \ref{lem:H_varphi} and \ref{lem:varphi'_terminates}. Given some $\Sigma^0_1$-formula $\varphi(n, x, Z)$, we construct $\varphi'(n, x)$, which is supposed to hold if and only if $\varphi(n, x, Z)$ does for $Z$ with $H_{\varphi}(\restr{X}{x}, Z)$. In $\varphi'(n, x)$, we search for the truth value of $\varphi(n, x, Z)$ by evaluating terms. The termination of these evaluations will later be guaranteed by the existence of a well-order that is produced by the fact that~$\alpha \mapsto \alpha^X$ preserves well-orders.

\begin{definition}\label{def:varphi'}
	Given a well-order $X$, we define terms inductively using a new symbol $P$ as follows:
	\begin{itemize}
		\item $0$, $1$,
		\item $P(m, y, \tilde{s})$ for $m \in \N$ and $y \in X$ where $\tilde{s}$ is a finite sequence of numbers where the last entry can also be a term that has been defined previously.
	\end{itemize}
	We call terms that are different from $0$ or $1$ \emph{proper} terms.
	Given $\Sigma^0_1$-formulas $\varphi(n, x, Z)$ and $\psi(n, x, Z)$, we define $\Sigma^0_1$-formulas $\varphi'(n, x)$ and $\psi'(n, x)$ that construct a sequence of these terms.

	Consider $\Sigma^0_1$-formulas $\varphi_0(n, x, s)$ and $\psi_0(n, x, s)$ that satisfy Kleene's normal form for $\varphi(n, x, Z)$ and $\psi(n, x, Z)$, respectively (cf.~\cite[Theorem~II.2.7]{Simpson}). This~means that $\varphi(n, x, Z)$ holds if and only if $\varphi_0(n, x, s)$ does for an initial segment~$s$~of~$Z$. Moreover, $\varphi_0(n, x, s)$ implies $\varphi_0(n, x, t)$ for any extension $t$ of $s$. Analogously, the same holds for $\psi(n, x, Z)$ and $\psi_0(n, x, s)$. Now, our sequence begins with the term $P(n, x, \langle \rangle)$ and subsequent members are constructed by case distinction on the previous member as follows:
	\begin{itemize}
		\item $0$ or $1$: The construction terminates, no new member is added.
		\item $P(m, y, \tilde{s})$:
		\begin{itemize}
			\item $\tilde{s}$ is a sequence of numbers and $\varphi_0(m, y, \tilde{s})$ or $\psi_0(m, y, \tilde{s})$ holds:
			
			The next member of our sequence is $1$ if $\varphi_0(m, y, \tilde{s})$ holds and $0$ otherwise.
			\item $\tilde{s}$ is a sequence of numbers but neither $\varphi_0(m, y, \tilde{s})$ nor $\psi_0(m, y, \tilde{s})$ does hold:
			
			The next member is $P(m, y, \tilde{s}')$ with $\tilde{s}' := \tilde{s} * \langle P(m', y', \langle\rangle) \rangle$ if $|\tilde{s}|$ is the code of a pair $(m', y')$ with $y' <_X y$. Otherwise, we have $\tilde{s}' := \tilde{s} * \langle 0 \rangle$.
			\item $\tilde{s}$ has a proper term as its last entry:
			
			In this case, we have $\tilde{s} = \tilde{s}' * \langle t\rangle$ for some sequence $\tilde{s}'$ of number-terms and a proper term $t$. Let $t'$ be the term that results from $t$ if it makes a step using this case distinction recursively. We continue with the term $P(m, y, \tilde{s}' * \langle t'\rangle)$.
		\end{itemize}
	\end{itemize}
	The formula $\varphi'(n, x)$ holds if and only if both this construction terminates and the last sequence member is equal to $1$. Similarly, $\psi'(n, x)$ holds if and only if both this constructions terminates and the last sequence member is equal to $0$.
\end{definition}

\begin{lemma}\label{lem:H_varphi}
	Let $X$ be some well-order and let both $\varphi'(n, x)$ and $\psi'(n, x)$ be defined like in Definition~\ref{def:varphi'} for some $\Sigma^0_1$-formulas $\varphi(n, x, Z)$ and $\psi(n, x, Z)$.
	We assume (\ref{eq:assumption_SETR}) and that $\varphi'(n, x) \leftrightarrow \lnot\psi'(n, x)$ holds for all $n \in \N$ and $x \in X$. Then, we have $H_{\varphi}(X, Y)$ for the set $Y$ that contains a pair $(x, n)$ if and only if $x \in X$ and $\varphi'(n, x)$ hold.
\end{lemma}

\begin{proof}
	The formula ``$x \in X \land \varphi'(n, x)$'' is $\Delta^0_1$ by assumption. Therefore, we can construct a set $Y$ that contains exactly those pairs $(x, n)$ with $x \in X$ and $\varphi'(n, x)$ using $\Delta^0_1$-comprehension. For the claim, we have to show that $\varphi'(n, x)$ holds if and only if $\varphi(n, x, Y^x)$ does, for any $n \in \N$ and $x \in X$.
	
	We begin by showing that $\varphi'(n, x)$ implies $\varphi(n, x, Y^x)$. Assume that $\varphi'(n, x)$ holds. Therefore, the sequence constructed by $\varphi'$ is finite and ends with a $1$. With this, the penultimate member is $P(n, x, s)$ for some finite $0$-$1$-sequence $s$ with $\varphi_0(n, x, s)$. By definition of the construction and $Y$, $s$ is an initial segment of~$Y^x$. This can be seen as follows: By construction, $s_i$ is zero for any $i < |s|$ that either does not code a pair or codes a pair $(y, m)$ with $y \nless_X x$. For any $i < |s|$ that does code a pair $(y, m)$ with $y <_X x$, we add a new term $P(m, y, \langle\rangle)$. The evaluation now proceeds exactly like that of $\varphi'(m, y)$. After all necessary steps, $P(m, y, \langle\rangle)$ will have been replaced by $1$ if $\varphi'(m, y)$ holds and by $0$ if $\psi'(m, y)$ holds. By assumption, these evaluations terminate and the result can be decided.	
	We conclude $\varphi(n, x, Y^x)$.
	
	Now, we show that $\varphi(n, x, Y^x)$ implies $\varphi'(n, x)$. We assume for contradiction that this does not hold for some $n \in \N$ and $x \in X$. In order to apply a least number principle, we use induction along $X$ on the $\Pi^0_1$-formula ${\forall n\ \big(\varphi(n, x, Y^x) \to (x, n) \in Y\big)}$. The required induction principle is implied by $\Sigma^0_0$-induction along $\omega \cdot X$ (cf.~\cite[Lemma 4.5]{Sommer}) and, therefore, available in $\RCAn$. Let $x$ be the smallest element~${x \in X}$ such that $\varphi(n, x, Y^x)$ holds but $\varphi'(n, x)$ does not for some $n \in \N$. Since~$x$ is the smallest such value and we have already shown that $\varphi'(m, y)$ always implies $\varphi(m, y, Y^y)$, we know that $H_{\varphi}(\restr{X}{x}, Y^x)$ holds. Therefore, the formulas $\varphi(n, x, Y^x)$ and $\lnot\psi(n, x, Y^x)$ are equivalent by assumption (\ref{eq:assumption_SETR}). We summarize that both $\psi(n, x, Y^x)$ and $\varphi'(n, x)$ do not hold. Since $\varphi'(n, x) \leftrightarrow \lnot\psi'(n, x)$ holds for all $n \in \N$ and $x \in X$, we know that its constructed sequence ends with a $0$ and that it has $P(n, x, s)$ as its penultimate term for a sequence $s$ such that, therefore, $\psi_0(n, x, s)$ must hold. Similar to above, we know that $s$ is an initial segment of $Y^x$. Therefore, we have $\psi(n, x, Y^x)$, which is a contradiction.
\end{proof}

\begin{lemma}\label{lem:varphi'_terminates}
	Let $X$ be some well-order and let both $\varphi'(n, x)$ and $\psi'(n, x)$ be defined like in Definition \ref{def:varphi'} for some $\Sigma^0_1$-formulas $\varphi(n, x, Z)$ and $\psi(n, x, Z)$.
	Assuming (\ref{eq:assumption_SETR}) and that $\alpha \mapsto \alpha^X$ preserves well-orders, we can show that $\varphi'(n, x) \leftrightarrow \lnot\psi'(n, x)$ holds for all $n \in \N$ and $x \in X$.
\end{lemma}

If we replace the assumption that $\alpha \mapsto \alpha^X$ preserves well-orders by the assumption that $\Pi^0_2$-induction along $X$ holds, the proof of this lemma becomes straightforward and can be used as an additional proof that $\Pi^0_2$-induction along $X$ implies $\SETR_X$ (cf.~Corollary~\ref{cor:induction_to_SETR_X} for the first one): We simply assume that $x \in X$ is the least element of $X$ such that $\varphi'(n, x) \leftrightarrow \lnot\psi'(n, x)$ does not hold for all $n \in \N$. Since this equivalence holds for all elements less than $x$, we know that the sequence of terms constructed by $\varphi'(m, y)$ and $\psi'(m, y)$ always terminates for all $m \in \N$ and~$y <_X x$. From this, we can easily show that the construction must also terminate for $\varphi'(n, x)$ and~$\psi'(n, x)$.

\begin{proof}[Proof of Lemma~\ref{lem:varphi'_terminates}]
	We assume that $X$ is not empty, otherwise the claim is clear. We use the fact that $\alpha \mapsto \alpha^X$ preserves well-orders in order to show that the search of $\varphi'(n, x)$ and $\psi'(n, x)$ always terminates for all $n \in \N$ and $x \in X$. For this, we will assign an element of $T^X$ to each member of the sequence that they produce. Let $T(n, x)$ be defined as the binary tree that includes a finite $0$-$1$-sequence $s$ if and only if neither $\varphi_0(n, x, s)$ nor $\psi_0(n, x, s)$ holds, for any $n \in \N$ and $x \in X$. Then, the partial order $T$ consists of $0_T$, $1_T$, and tuples $(n, x, t, b)$ for $n \in \N$, $x \in X$, $t \in T(n, x)$, and $b \in 2$. The elements $0_T <_T 1_T$ are the smallest ones. The tuples are ordered by priority using the order of $\N$, $X$, $T(n, x)$, and finally $2$. In general, $T$ is not a linear order. Also, if there is some sequence $Z$ with neither $\varphi(n, x, Z)$ nor $\psi(n, x, Z)$ for some $n \in \N$ and $x \in X$, it is not well-founded. Later, we will use a trick to circumvent both of these problems.
	
	For proof by contradiction, let $(f_i)_{i \in \N}$ be the infinite sequence of terms that $\varphi'(n, x)$ generates. Now, we are ready to define a similar sequence $(g_i)_{i \in \N}$ of elements in $T^X$. We define the element $\beta(t)$ of $T^X$ depending on the term $t$ by the following case distinction:
	\begin{itemize}
		\item $t \in \{0, 1\}$: Here, we set $\beta(t) := 0_{T^X}$ (i.e.~$\langle \rangle \in T^X$).
		\item $t = P(m, y, \tilde{s})$:
		\begin{itemize}
			\item $\tilde{s}$ is a sequence of numbers and $\varphi_0(m, y, \tilde{s})$ or $\psi_0(m, y, \tilde{s})$ holds:
			
			We define $\beta(t) := 1$ (i.e.~$\langle (0_X, 1_T) \rangle \in T^X$ where $0_X$ is the smallest element of $X$).
			\item $\tilde{s}$ is a sequence of numbers but neither $\varphi_0(m, y, \tilde{s})$ nor $\psi_0(m, y, \tilde{s})$ does hold:
			
			We define $\beta(t) := T^y \cdot (m, y, \tilde{s}, 1_2)$ (i.e.~$\langle (y, (m, y, \tilde{s}, 1_2)) \rangle \in T^X$).
			\item $\tilde{s}$ has a proper term as its last entry:
			
			Let $\tilde{s} = \tilde{s}' * \langle t'\rangle$ where $\tilde{s}'$ is a sequence of number-terms and $t'$ is a proper term. We define $\beta(t) := T^y \cdot (m, y, \tilde{s}', 0_2) + \beta(t')$\\(i.e.~${\langle (y, (m, y, \tilde{s}', 0_2))\rangle * \beta(t') \in T^X}$).
		\end{itemize}
	\end{itemize}
	The sequence $\beta$ that is defined as concatenation in the last case must be an element of $T^X$ since any term $t' = P(m', y', \tilde{s}'')$ that might appear as the last element of $\tilde{s}$ satisfies $y' <_X y$. Therefore, we know that for any pair $(y', t) \in \beta(t')$ the inequality~$y' <_X y$ must hold.
	
	Using this mapping from terms to $T^X$, we can define the sequence $(g_i)_{i \in \N}$ such that each $g_i$ is defined as $\beta(f_i) \in T^X$ for each $i \in \N$. Comparing the steps of our term evaluation with these assignments to elements in~$T^X$, one can verify that $(g_i)_{i \in \N}$ is a descending sequence in $T^X$.
	
	As already mentioned, we cannot use the assumption that the mapping $\alpha \mapsto \alpha^X$ preserves well-orders at this moment, since it is possible that $T$ is not a well-order. For this, we use the following trick: We define a new partial order $T'$ consisting of pairs $(i, t)$ for $i \in \N$ and $t \in T$ where $i$ is the smallest index such that $t$ occurs in~$g_i$, i.e., $(x, t) \in g_i$ holds for some $x \in X$. The order of two elements $(i, t)$ and~$(i', t')$ is simply given by that of $t$ and $t'$ with respect to $T$. Trivially, we know that~$(g_i)_{i \in \N}$ also lives in $T'^X$ (or that a very similar sequence $(g'_i)_{i \in \N}$ does so where each pair~$(x, t)$ in $g_i$ is replaced by $(x, (j, t))$ for the unique $j \in \N$ with $(j, t) \in T'$).
	
	Now, we show that $T'$ is a linear order. For this, consider two incomparable elements $(i, (m, y, s, b))$ and $(i', (m', y', s', b'))$. We immediately know $m = m'$ and~$y = y'$. Otherwise, they would easily be comparable. Assume that $s$ and $s'$ are incomparable sequences in the tree $T(m, y)$. This implies that there is some index~$k < \min(|s|, |s'|)$ with $s_k \neq s'_k$. By construction of $T'$, we know that both terms~$(m, y, s, b)$ and $(m, y, s', b')$ appear somewhere in the sequence $(g_i)_{i \in \N}$. By definition of this sequence, we can find both $P(m, y, \tilde{s})$ and $P(m, y, \tilde{s}')$ as subterms in $(f_i)_{i \in \N}$ where $\tilde{s}$ and $\tilde{s}'$ are defined like $s$ and $s'$ but may also additionally include a proper term as last member. Since one of $\tilde{s}_k$ or $\tilde{s}'_k$ is equal to one, we know that~$k$ is the code of a pair $(l, z)$ with $z <_X y$. Moreover, we know that during the computation of $\tilde{s}$, the term $P(l, z, \langle \rangle)$ was evaluated to $\tilde{s}_k$. Similarly, during the computation of $\tilde{s}'$, the term $P(l, z, \langle \rangle)$ was evaluated to~$\tilde{s}'_k$. However, any evaluation of a term (as long as it terminates) must produce the same value, which can be seen from the definition of our sequence $(f_i)_{i \in \N}$. This is a contradiction, and we conclude $s = s'$. Finally, we know that $b$ and $b'$ cannot be incomparable. We have shown that $T'$ is a linear order.
	
	With the infinite descending sequence $(g_i)_{i \in \N}$, we have already seen that $T'^X$ cannot be a well-order. Now, we apply the principle that $\alpha \mapsto \alpha^X$ preserves well-orders. With this, we know that the linear order $T'$ cannot be a well-order, either. Consider an infinite descending sequence $(h_i)_{i \in \N}$ in $T'$. We identify this with the uniquely corresponding infinite descending sequence $(m_i, y_i, s^i, b_i)_{i \in \N}$ in $T$ (we use superscripts $s^i$ in order to avoid confusion because ``$s_i$'' could also be understood as the $i$-th member of some sequence $s$). Since $\N$ is a well-order, we know that at some point, the sequence $(m_i)_{i \in \N}$ cannot change its value anymore. After this, the same holds for the well-order $X$ and the sequence $(y_i)_{i \in \N}$. We, therefore, simply assume w.l.o.g.~that our sequence starts from that point and, furthermore, that all members of both sequences $(m_i)_{i \in \N}$ and $(y_i)_{i \in \N}$ are equal to fixed values $m \in \N$ and $y \in X$, respectively. We know that $(b_i)_{i \in \N}$ is infinitely often equal to $0$ or to~$1$. Similarly to above, we w.l.o.g.~continue with a subsequence such that $(b_i)_{i \in \N}$ is constant. Now that $(m_i)_{i \in \N}$, $(y_i)_{i \in \N}$, and $(b_i)_{i \in \N}$ are constant sequences, we know that $(s^i)_{i \in \N}$ must be an infinite descending chain in $T(m, y)$.
	
	We have an infinite descending chain $(s^i)_{i \in \N}$ of sequences in $T(m, y)$. Since we started with an infinite descending sequence in $T'$, we know that $(f_i)_{i \in \N}$ must contain terms $P(m, y, \tilde{s}^i)$ for each index $i \in \N$ (where $\tilde{s}^i$ is defined like $s^i$ but may have an additional proper term as last member). The existence of these terms in~$(f_i)_{i \in \N}$ implies that for any index $k < |s^i|$ such that $k$ is the code of a pair~$(m', y')$ with~$y' <_X y$ we know that the evaluation of $P(m', y', \langle \rangle)$ terminates. Since~$(s^i)_{i \in \N}$ is an infinite descending sequence (ordered by extension), its elements~$s^i$ get arbitrarily long. We conclude that the evaluation of $P(m', y', \langle \rangle)$ terminates for all pairs $(m', y')$ with $y' <_X y$. Thus, we have $\varphi'(m', y') \leftrightarrow \lnot\psi'(m', y')$ for any $m' \in \N$ and $y' <_X y$.
	
	Using Lemma \ref{lem:H_varphi}, we know that $H_{\varphi}(\restr{X}{y}, Y)$ holds for the set $Y$ that contains a pair $(m', y')$ if and only if $y' \in \restr{X}{y}$, i.e.~$y' <_X y$, and $\varphi'(m', y')$ hold. With $H_{\varphi}(\restr{X}{y}, Y)$, however, we know by assumption (\ref{eq:assumption_SETR}) that $\varphi(m, y, Y)$ holds if and only if $\lnot\psi(m, y, Y)$ does. We conclude that the sequence $(f_i)_{i \in \N}$ cannot contain all terms $P(m, y, \tilde{s}^i)$ for our sequence $(s^i)_{i \in \N}$. Otherwise, this would mean that neither $\varphi_0(m, y, s^i)$ nor $\psi_0(m, y, s^i)$ does hold for any index $i \in \N$. By definition of $\varphi_0$ and $\psi_0$, however, this would imply that neither $\varphi(n, x, Y)$ nor $\psi(n, x, Y)$ does hold. Recall for this, that any $s^i$ is an initial segment of $Y$ by construction. Since one of $\varphi(n, x, Y)$ or $\psi(n, x, Y)$ must hold under the assumption $H_{\varphi}(\restr{X}{y}, Y)$, we have arrived at a contradiction.
\end{proof}

\begin{proof}[Proof of Proposition \ref{prop:alpha_x_proves_SETR_X}]
	Assume that $\alpha \mapsto \alpha^{X}$ preserves well-orders. We apply Lemma \ref{lem:varphi'_terminates} and know that $\varphi'(n, x) \leftrightarrow \lnot\psi'(n, x)$ holds for all $n \in \N$ and $x \in X$. With $\Delta^0_1$-comprehension, we construct a set $Y$ that contains $(x, n)$ if and only if~$x \in X$ and $\varphi'(n, x)$ hold. With Lemma \ref{lem:H_varphi}, we conclude $H_{\varphi}(X, Y)$.
\end{proof}

\section{Weak Effective Transfinite Recursion}\label{sec:WETR}

In this section, we define Weak Effective Transfinite Recursion as considered in~\cite{Dzhafarov}. Then, in Proposition \ref{prop:WKL_proves_WETR'}, we show that $\WKL_0$ proves a variation of the principle where we even allow well-founded relations instead of only well-orders. After that, in Proposition \ref{prop:RCA_proves_WETR_rule}, we demonstrate that $\RCAn$ can prove a rule version of $\WETR$. Then, in Proposition \ref{prop:WETR_proves_WKL}, we show that $\WETR_X$, i.e.~$\WETR$ for a certain well-order $X$, implies the disjunction of $\WKL$ and $\SETR_X$. Together with the previous results, this yields an equivalence.

\begin{definition}[Weak Effective Transfinite Recursion]\label{def:WETR}
	Consider some well-order~$X$ together with $\Sigma^0_1$-formulas $\varphi(n, x, Z)$ and $\psi(n, x, Z)$. The principle $\WETR_X$ of \emph{Weak Effective Transfinite Recursion} along $X$ is defined as follows: If we have
	\begin{equation}\label{eq:delta_0_1}
		\forall n \in \N\ \forall x \in X\ \forall Z \subseteq \N\ \big(\varphi(n, x, Z) \leftrightarrow \lnot \psi(n, x, Z)\big)\komma
	\end{equation}
	then there exists some set $Y$ with $H_{\varphi}(X, Y)$. We write $\WETR$ for the principle that $\WETR_X$ holds for any well-order $X$.
\end{definition}

Note that Definition \ref{def:WETR} is similar to Definition \ref{def:SETR}. The only difference is that (\ref{eq:delta_0_1}) is missing the premise $H_\varphi(\restr{X}{x}, Z)$ in contrast to (\ref{eq:assumption_SETR}). Similar to our discussion in Remark \ref{rem:single_formula} for $\SETR_X$ and $\SETR$, we can write the principles $\WETR_X$ and $\WETR$ as single formulas.
The main result of this section is the following theorem:

\begin{theorem}[$\RCAn$]\label{thm:WETR}
	Consider some well-order $X$. Then, the following are equivalent:
	\begin{enumerate}[label=\alph*)]
		\item $\WETR_X$
		\item The disjunction of $\WKL$ and $\SETR_X$.
	\end{enumerate}
\end{theorem}

We begin by showing that $\WKLn$ proves $\WETR$. In fact, using $\WKL$, we can prove the (seemingly) stronger principle of $\WETR$ \emph{for well-founded relations}:

\begin{definition}[$\WETR$ for well-founded relations]
	We define \emph{$\WETR$ for well-founded relations} ($\WETR'$) like $\WETR$ except that $X$ may now be any well-founded relation, i.e., $X$ no longer has to be a linear order.
\end{definition}

\begin{proposition}\label{prop:WKL_proves_WETR'}
	$\WKLn$ proves $\WETR'$.
\end{proposition}

\begin{proof}
	Consider an arbitrary well-founded relation $X$ with $\Sigma^0_1$-formulas $\varphi(n, x, Z)$ and $\psi(n, x, Z)$ that form an instance of $\WETR'$, i.e., we assume (\ref{eq:delta_0_1}).
	We can find $\Delta^0_0$-formulas $\varphi_0$ and $\psi_0$ satisfying Kleene's normal form for $\varphi$ and $\psi$, respectively, similar to the construction in Definition \ref{def:varphi'}: $\varphi(n, x, Z)$ holds if and only if $\varphi_0(n, x, s)$ does for an initial segment $s$ of $Z$. Moreover, $\varphi_0(n, x, s)$ implies $\varphi_0(n, x, t)$ for any extension $t$ of $s$. The properties of $\psi_0$ are analogous.
	
	Also, recall the definition of $T(n, x)$: For every $n \in \N$ and $x \in X$ we define $T(n, x)$ to be the tree consisting of all finite $0$-$1$-sequences $s$ such that neither $\varphi_0(n, x, s)$ nor $\psi_0(n, x, s)$ holds. In $\WKLn$, each tree $T(n, x)$ must be finite. Otherwise, assume that $T(n, x)$ is infinite for some $n$~and~$x$. Then, it must contain an infinite path~$Z$. We conclude that neither $\varphi(n, x, Z)$ nor $\psi(n, x, Z)$ does hold since neither $\varphi_0(n, x, Z[m])$ nor $\psi_0(n, x, Z[m])$ does hold for any $m$, which is a contradiction to the assumptions of $\WETR'$.
	We define a partial function $t$ that gets a pair $(n, x)$ with $n \in \N$ and $x \in X$ and searches for an upper bound on the length of sequences in $T(n, x)$, i.e., we have $t(n, x) >_\N |s|$ for any $s \in T(n, x)$ for $n \in \N$ and $x \in X$. The underlying set of $t$ can be constructed in $\RCAn$ since it is decidable whether some number is such a bound. Moreover, $t$ is total since all trees $T(n, x)$ are finite. With this, we know that for any sequence $s$ with $|s| \geq t(n, x)$ either $\varphi_0(n, x, s)$ or $\psi_0(n, x, s)$ must hold.
	
	Next, we define a program $e$ that takes two elements $n \in \N$ and $x \in X$ and computes the following: We wait for the termination of $e(m, y)$ for all pairs $(m, y) <_{\N} t(n, x)$  with $y <_X x$, where we interpret $(m, y)$ as a natural number. After this, we collect the results in a $0$-$1$-sequence $s$ of length $|s| = t(n, x)$:
	\begin{equation}\label{eq:sequence_s}
		s_i :=
		\begin{cases}
			e(m, y) & \text{ if } i = (m, y) \text{ for } y <_X x \text{,}\\
			0 & \text{ otherwise.}
		\end{cases}
	\end{equation}
	After this construction, $e(n, x)$ terminates. It assumes the value $1$ if $\varphi_0(n, x, s)$ is true. Otherwise, its value is $0$.
	
	In order to show that $e$ is a total program with respect to its domain $\N \times X$, i.e., that this construction terminates for $(n, x)$ (and all instances $(m, y)$ that are used by the recursion), we consider a tree $T$ of pairs in $\N \times X$ (interpreted as natural numbers). The tree is constructed as follows: It contains the root $\langle\rangle$ and every finite sequence $s$ of elements in $\N \times X$ starting with $(n, x)$ such that for each index $0 \leq i < |s| - 1$, and pairs $(m, y) := s_i$ and $(m', y') := s_{i+1}$, we have both $y' <_X y$ and $(m', y') <_{\N} t(m, y)$. Intuitively, each node $(m, y)$ in our tree has all pairs $(m', y')$ as children for which $e(m, y)$ uses the value of $e(m', y')$ during its computation.
	
	In these next steps, we want to prepare everything for an application of the \emph{Bounded} K\H{o}nig's Lemma. A tree is bounded if and only if there exists some bounding function $g: \N \to \N$ such that for each sequence $s$ in our tree, we have $s_i \leq g(i)$ for all $i < |s|$ (cf.~\cite[Definition IV.1.3]{Simpson}). Clearly, it is decidable whether a sequence lies in our tree $T$ or not. Moreover, the tree is bounded by a function using the following construction: There is only one sequence of length $0$ and only one of length $1$ (trivially bounded by $(n, x)$ itself). For length $k > 1$ consider the bound $b$ that we computed for sequences of length $k - 1$. This means that all sequences $s$ of length $k - 1$ end with a pair $(m, y) \leq_{\N} b$. By construction, we know that all children of the node represented by $s$ are pairs $(m', y')$ with $(m', y') <_{\N} t(m, y)$. Therefore, all children of $s$ are bounded by $t(m, y)$. We conclude that the maximum of $t(m, y)$ for all pairs $(m, y) \leq_\N b$ is a bound for length $k$.
	
	Now, we have a bounded tree. Next, we show that $T$ has no infinite path. Otherwise, this path would be a sequence of pairs $((n_i, x_i))_{i \in \N}$. By construction of our tree, we have $x_{i+1} <_X x_i$ for all $i \in \N$, which contradicts our assumption that $X$ is well-founded. With this, we can apply Bounded K\H{o}nig's Lemma, which is implied by Weak K\H{o}nig's Lemma (cf.~\cite[Lemma IV.1.4]{Simpson}).
	
	We now know that $T$ is finite. Assume that $e(m, y)$ is not defined for some pair $(m, y)$ that occurs as last element of some sequence $s$ (not necessarily a leaf) in our tree. Since the tree is finite, we can assume that $s$ has the greatest length of such a counterexample using $\Sigma^0_1$-induction. This implies that $e$ is defined for all pairs $(m', y')$ with $(m', y') <_{\N} t(m, y)$ and $y' <_X y$. But then, we can compute the value of $e(m, y)$. Therefore, $e$ is defined for all pairs $(m, y)$ occurring in our tree, including $(n, x)$.
	
	We can now produce some $Y$ consisting of all pairs $(x, n)$ with $x \in X$ and $e(n, x) = 1$ using $\Delta^0_1$-comprehension. In the final steps of this proof, we show that $H_{\varphi}(X, Y)$ holds for this $Y$. In order to do this, we simply have to make sure that $e(n, x) = 1$ holds if and only if $\varphi(n, x, Y^x)$ does.
	
	If $e(n, x) = 1$ holds, then $\varphi_0(n, x, s)$ holds for $s$ defined as in (\ref{eq:sequence_s}). By definition of $\varphi_0$ and since $s = Y^x[t(n, x)]$ holds, we have $\varphi(n, x, Y^x)$.
	Conversely, if $e(n, x) = 1$ does not hold, then we have $e(n, x) = 0$ since $e$ is total on its domain~${\N \times X}$. With this, $\varphi_0(n, x, s)$ does not hold. By definition of $t(n, x)$, we know that $\psi_0(n, x, s)$ must hold since we have $|s| = t(n, x)$. Therefore, we have $\psi(n, x, Y^x)$ and finally~$\lnot \varphi(n, x, Y^x)$.
\end{proof}

This result immediately yields one direction of our theorem:

\begin{corollary}[$\RCAn$]
	Consider some well-order $X$. The disjunction of $\WKL$ and $\SETR_X$ implies $\WETR_X$.
\end{corollary}

\begin{proof}
	Assume that we have $\WKL$. Then, by Proposition \ref{prop:WKL_proves_WETR'}, this implies $\WETR'$, $\WETR$, and finally $\WETR_X$. Assume that we have $\SETR_X$. We immediately have $\WETR_X$ since $\SETR_X$ has weaker assumptions on $\varphi(n, x, Z)$ than $\WETR_X$.
\end{proof}

We continue with a rule that is admissible over $\RCAn$: We present a variant of $\WETR$ where we assume that the decidability of $\varphi$ is provable (and not just an assumption). With this (and restricting ourselves to linear orders $X$), we can eliminate both occurrences of Weak K\H{o}nig's Lemma from the proof of Proposition~\ref{prop:WKL_proves_WETR'}. This yields a rule that is admissible over $\RCAn$:

\begin{proposition}[Rule Version of $\WETR$]\label{prop:RCA_proves_WETR_rule}
	Let $P$ be an arbitrary $\Sigma^1_1$-formula, and let both $\varphi(n, x, Z)$ and $\psi(n, x, Z)$ be $\Sigma^0_1$-formulas . If we can derive
	\begin{align}\label{eq:RCAn_proves_P_implies_totality}
		&\RCAn \vdash P \to \forall n \in \N\ \forall x \in X\ \forall Z \subseteq \N\ \big(\varphi(n, x, Z) \leftrightarrow \lnot \psi(n, x, Z)\big)\komma\\
	\intertext{then this implies}
		&\RCAn \vdash P \land \WO(X) \to \exists Y\ H_\varphi(X, Y)\period\notag
	\end{align}
\end{proposition}

\begin{proof}
	Note that we used two instances of Weak K\H{o}nig's Lemma in our proof of Proposition \ref{prop:WKL_proves_WETR'}: First, we used it to show that all trees $T(n, x)$ are finite, and then we applied it in order to derive that the tree constructed for the termination of $e(n, x)$ is also finite. Now, the idea of this proof is to eliminate both these instances of $\WKL$.
	
	As we have seen in the proof of Proposition \ref{prop:WKL_proves_WETR'}, we can show in $\WKLn$ that the function $t(n, x)$ finishes its search for each $n \in \N$ and $x \in X$ if $\varphi(n, x, Z) \leftrightarrow \lnot \psi(n, x, Z)$ holds for any $n$, $x$, and $Z$. Note that for this we did not use the assumption that $X$ is well-founded. Moreover, the definition of $t(n, x)$ only relies on the definitions of $\varphi_0$ and $\psi_0$, and no further information that might only be available in the presence of $\WKL$. Combining this with our assumption (\ref{eq:RCAn_proves_P_implies_totality}), we know that $\WKLn$ proves the totality of $t$ under the assumption $P$. Totality is a $\Pi^0_2$-statement. Thus, the implication from $P$ to the totality of $t$ is a $\Pi^1_1$-statement. We can, therefore, since $\WKLn$ is $\Pi^1_1$-conservative over $\RCAn$ (cf.~\cite[Corollary~IX.2.6]{Simpson}) already prove this implication in $\RCAn$.
	
	Next, we want to show without the use of $\WKL$ that the tree used for the termination of $e(n, x)$ is already finite. This is possible since we are now working with well-orders, where we can find the maximum of finite non-empty subsets of~$X$, instead of just well-founded relations. Assume for contradiction that this tree is infinite. Since the tree is bounded, we can for each $k > 0$ compute the maximum of all $y$ in pairs $(m, y)$ that are the last element of sequences of length $k$ in our tree. This produces a sequence $(y_k)_{k > 0}$. Moreover, this sequence is descending: For each member $y_{k+1}$ with $k > 0$, we know that there exists some sequence $s$ of length $|s| = k + 1$ in our tree that ends with the pair $(m, y_{k+1})$ for some $m \in \N$. The penultimate pair in $s$ must, therefore, be some pair $(m', y)$ with $m \in \N$ and~${y > y_{k+1}}$. We conclude $y_k \geq y > y_{k+1}$. Now, we have produced an infinite descending sequence in $X$, which contradicts our assumption that $X$ was well-founded.
	
	Finally, we conclude that both occurrences of Weak K\H{o}nig's Lemma can be replaced: The first one by the assumption $P$ and the second one by the assumption that $X$ is a linear order. Using the rest of the proof for Proposition \ref{prop:WKL_proves_WETR'} that only relies on $\RCAn$, we have shown our claim.
\end{proof}

\begin{remark}
	We show that the rule above cannot be extended to relations $X$ that are only well-founded partial orders. Consider an arbitrary function $f: \N \to \N$ and the following formula $\varphi(n, x, Z)$:
	\begin{equation*}
		\varphi(n, x, Z) := \exists m \leq_\N x\ f(m) = n \lor (n, x * \langle 0 \rangle) \in Z \lor (n, x * \langle 1 \rangle) \in Z\period
	\end{equation*}
	We see that this formula lives in $\Delta^0_1$ and our rule from Proposition \ref{prop:RCA_proves_WETR_rule} tells us that for any well-order $X$, there exists some set $Y$ with $H_{\varphi}(X, Y)$.
	
	Now, assume that this rule does not only work for well-orders $X$ but also for well-founded partial orders. Using this assumption, we want to show $\WKL$. In order to derive a contradiction, assume $\lnot \WKL$. Now, we can use the infinite tree~$T$ without path as our well-founded relation by defining $s >_T t$ whenever $t$ is a proper extension of $s$ for $s, t \in T$.
	
	For any $n \in \N$, we want to show that $(n, \langle\rangle) \in Y$ holds if and only if $n$ is in the range of $f$. Assume that $(n, s) \in Y$ holds for some sequence $s \in T$ but $n$ is not in the range of $f$. By definition of $\varphi(n, x, Z)$, one of $(n, s * \langle 0 \rangle) \in Y$ or $(n, s * \langle 1 \rangle) \in Y$ must hold. Thus, we can construct a path in $T$. This contradiction implies that $n$ must lie in the range of $f$. For the other direction, assume that $n$ is in the range of $f$ witnessed by $m \in \N$. Since $T$ is infinite, there must be some sequence $s \in T$ with a code greater than the number $m$. We conclude $\varphi(n, s, Y^s)$ and, therefore, $(n, s) \in Y$. Using induction on the definition of $\varphi(n, x, Z)$, we see that $(n, s') \in Y$ holds for any initial segment $s'$ of $s$. This includes $s' := \langle \rangle$.

	Now, using $\Delta^0_1$-comprehension, we can produce the ranges of arbitrary functions~$f: \N \to \N$. This leads to $\ACAn$ (cf.~\cite[Lemma~III.1.3]{Simpson}), which is a contradiction since it implies $\WKLn$. We conclude that the existence of a set $Y$ with $H_{\varphi}(T, Y)$ implies $\WKL$.
\end{remark}

We finish the proof of Theorem \ref{thm:WETR} with the following proposition:

\begin{proposition}[$\RCAn$]\label{prop:WETR_proves_WKL}
	$\WETR_X$ implies $\WKL$ or $\SETR_X$.
\end{proposition}

\begin{proof}
	We derive $\Pi^0_2$-induction along $X$ using $\WETR_X$ and $\lnot \WKL$. Using Corollary \ref{cor:induction_to_SETR_X} or both Propositions \ref{prop:induction_proves_alpha_x} and \ref{prop:alpha_x_proves_SETR_X}, this will prove our claim. The latter assumption allows us to code witnesses for the induction into the set that is constructed by the recursion principle. In general, this is not possible but in the context of $\lnot\WKL$, there exists a continuous surjective function from $2^\N$ to $\N$, which enables us to do just that. The proof is very similar to that of Proposition \ref{prop:SETR_X_implies_induction}, where we did not need explicit witnesses but simply ensured the existence of such witnesses using the assumption $H_\varphi(\restr{X}{x}, Z)$ at each stage.
	
	We assume that the $\Pi^0_2$-formula $\forall n\ \exists m\ P(x, n, m)$ is progressive in $x$ where $P(x, n, m)$ is some $\Delta^0_0$-formula. Under this assumption, we want to show that it holds true for any $x \in X$. We require that for any $x \in X$ and $n \in \N$, there is at most one $m \in \N$ with $P(x, n, m)$. This can easily be achieved by replacing $P(x, n, m)$ by the predicate ${P(x, n, m) \land \forall m' <_\N m\ \lnot P(x, n, m')}$.

	Let $T$ be a binary tree that is both infinite and has no path. Using $\Delta^0_1$-comprehension, we can construct the set $L \subseteq \N$ of lengths $l \in \N$ such that there exists a sequence $s \notin T$ with $|s| = l$ where all proper initial segments of $s$ are in $T$. This set $L$ is infinite: For any $n \in \N$, we can find a sequence $s \in T$ with $|s| \geq_\N n$ since $T$ is infinite. Now, we extend $s$ with zeros until we have reached an extension $s' \in T$ such that $s' * \langle 0 \rangle \notin T$ holds. This will happen eventually since $T$ has no paths. We conclude that we have found a sequence $s' \in T$ with $|s'| \geq_\N n$ and $|s'| + 1 \in L$ for an arbitrary $n \in \N$. The set $L$ is, therefore, infinite. Finally, we can find a surjective function $f: L \to \N$.
	
	Now, we can map arbitrary infinite $0$-$1$-sequences to natural numbers in the following way: Search for the shortest initial segment $s \notin T$ of the infinite sequence. This search will terminate as $T$ contains no paths. The natural number is then given by $f(|s|)$. Conversely, by definition of $L$ and the surjectivity of $f$, we can find such an (infinite) $0$-$1$-sequence for any $n \in \N$.
	
	Let $seq: \N \to \{0, 1\}^*$ be the function that maps any $n \in \N$ to the sequence $s \notin T$ with all proper initial segments in $T$ that has the smallest code such that $f(|s|) = n$ holds. From this point on, we also define $s_i := 0$ for any sequence $s \in \{0, 1\}^*$ and index $i \geq_\N |s|$.
	Furthermore, let $init: \{0, 1\}^* \to \{0, 1\}^*$ be the function that maps any sequence $s \neq \langle\rangle$ to its longest proper initial segment, i.e., we have $init(s) * \langle b \rangle = s$ for some $b \in \{0, 1\}$. The value of $init$ for the empty sequence $\langle \rangle$ is not important, we may simply set it to $\langle \rangle$. Now, we continue with the following definitions:
	\begin{align*}
		Z_{x, n} &:= \ \{i \in \N \mid (x, (n, i)) \in Z\}\komma\\
		Z \mapsto m &:\Longleftrightarrow \ \exists s \in \{0, 1\}^*\ s = Z[|s|] \land{} s \notin T \land init(s) \in T \land f(|s|) = m\period
	\end{align*}
	Later, we will use the combined expression $Z_{x, n} \mapsto m$ for sets $Z$ and $x, n, m \in \N$. Intuitively, this formula holds if and only if the $n$-th number coded into stage $x$ of~$Z$ is equal to $m$.
	
	We show that $Z \mapsto m$ is a $\Delta^0_1$-formula and that for any set $Z$ there is a unique~$m$ such that this formula holds. For the first claim, we simply prove that there is a unique sequence $s \in \{0, 1\}^*$ with $s = Z[|s|]$, $s \notin T$, and $init(s)$. The existence of this sequence is clear: Since $T$ has no infinite path, there must be some initial segment of $Z$ that is not included in $T$. We simply choose $s$ to be the shortest such sequence. For the uniqueness, consider a second sequence $t \in \{0, 1\}^*$ with these properties. As both $s$ and $t$ are initial segments of $Z$, one of them has to be a proper extension of the other. If $s$ is a proper initial segment of $t$, then $s \notin T$ and $init(t) \in T$ contradict each other. If $t$ is a proper initial segment of $s$, then $t \notin T$ and $init(s) \in T$ contradict each other. For the second claim, we simply consider the unique sequence $s \in \{0, 1\}^*$ from above and the value $m$ that $f(|s|)$ maps to. Should there be another possible value $m' \in \N$ with $Z \mapsto m'$ then that would correspond to a sequence $t \in \{0, 1\}^*$ with $f(|t|) = m'$ such that $t$ is an initial segment of $Z$ with $t \notin T$ and $init(t) \in T$. Using the arguments from before, this cannot be.
	
	Using the knowledge that for any set $Z$ there is a unique number $m \in \N$ with $Z \mapsto m$, we will write $P(x, n, F(Z))$ as abbreviation for the $\Delta^0_1$-statement\footnote{This $\Sigma^0_1$-formula is equivalent to the $\Pi^0_1$-formula ${\forall m\ (Z \mapsto m \to P(x, n, m))}$.} ${\exists m\ (Z \mapsto m \land P(x, n, m))}$ with $x \in X$ and $n \in \N$. Intuitively, $F$ can be interpreted as a function from sets to numbers that maps any set $Z$ to the unique $m$ with $Z \mapsto m$.
	
	Now, we are ready to define our $\Sigma^0_1$-formulas $\varphi((n, i), x, Z)$ and $\psi((n, i), x, Z)$. Here, we directly use pairs $(n, i)$ as our parameter:
	\begin{align*}
		\varphi((n, i), x, Z) := \ &\exists m\ \big(P(x, n, m) \land seq(m)_i = 1 \land{}\\
		&\ \forall (y, n') \leq_\N m\ \big(y <_X x \to P(y, n', F(Z_{y, n'}))\big)\big)\\
		\lnot\psi((n, i), x, Z) := &\ \forall m\ (P(x, n, m) \to seq(m)_i = 1) \land{}\\
		&\ \forall y <_X x\ \forall n'\ \big(\forall m <_\N (y, n')\ \lnot P(x, n, m) \to P(y, n', F(Z_{y, n'}))\big)
	\end{align*}
	These definitions are quite similar to those from the proof of Proposition \ref{prop:SETR_X_implies_induction}. Consequently, the proof of $\varphi((n, i), x, Z) \leftrightarrow \lnot\psi((n, i), x, Z)$ works analogous. However, instead of proceeding by case distinction on the truth of $\varphi((n, i), x, Z)$, we use case distinction with respect to
	\begin{equation}\label{eq:varphi_without_R}
		\exists m\ \big(P(x, n, m) \land \forall(y, n') \leq_\N m\ \big(y <_X x \to P(y, n', F(Z_{y, n'}))\big)\big)\komma
	\end{equation}
	i.e., $\varphi((n, i), x, Z)$ but without the conjunct $seq(m)_i = 1$.
	
	First, assume that (\ref{eq:varphi_without_R}) holds. If $\varphi((n, i), x, Z)$ does not hold, then there is some $m \in \N$ with $P(x, n, m)$ but not $seq(m)_i = 1$. We conclude that $\lnot \psi((n, i), x, Z)$ cannot hold. Now, assume that $\varphi((n, i), x, Z)$ does hold. We conclude ${\forall m\ \big(P(x, n, m) \to seq(m)_i = 1\big)}$ since there is at most one $m \in \N$ with $P(x, n, m)$. Let $m \in \N$ be a witness of $\varphi((n, i), x, Z)$. Under the assumption $\forall m <_\N (y, n')\ \lnot P(x, n, m)$, we have $(y, n') \leq_\N m$. Using $y <_X x$, we conclude $P(y, n', F(Z_{y, n'}))$.
	
	Secondly, assume that (\ref{eq:varphi_without_R}) does not hold. Here, we immediately see that $\varphi((n, i), x, Z)$ cannot hold, and we want to show that $\lnot \psi((n, i), x, Z)$ does not hold, as well. Assume that we have $P(y, n', F(Z_{y, n'}))$ for all $y <_X x$ and $n'$. By progressivity, this implies the existence of some $m$ with~$P(x, n, m)$. We conclude that (\ref{eq:varphi_without_R}) must hold, which is a contradiction. Thus, there is some pair $(y, n')$ with $y <_X x$ such that we do not satisfy $P(y, n', F(Z_{y, n'}))$. Using $\Pi^0_1$-induction, we can assume that $(y, n')$ has the smallest code with this property. By assumption $\lnot$(\ref{eq:varphi_without_R}), we now know that any $m \in \N$ with $P(x, n, m)$ must be greater or equal to the code of $(y, n')$. Thus, the premise ${\forall m <_\N (y, n')\ \lnot P(x, n, m)}$ in $\lnot \psi((n, i), x, Z)$ is satisfied. Since $P(y, n', F(Z_{y, n'}))$ does not hold by assumption on $(y, n')$, we conclude $\psi((n, i), x, Z)$.
	
	Now, we can apply $\WETR_X$, which produces a set $Y$ with $H_\varphi(X, Y)$. Consider the function $g: X \times \N \to \N$ with $((x, n), m) \in g$ if and only if $Y_{x, n} \mapsto m$. The fact that this is a $\Delta^0_1$-formula makes $g$ definable using $\Delta^0_1$-comprehension. The existence and uniqueness of $m$ ensures that $g$ is a function. Assume that $g$ is not a witness for the result of the induction, i.e., there is some pair $(x, n) \in X \times \N$ such that $P(x, n, g(x, n))$ does not hold. We assume that $(x, n)$ is the smallest such pair with respect to the lexicographical order on $X \times \N$. We know that $P(y, k, g(y, k))$ holds for all $y <_X x$ and $k \in \N$. By progressivity, it follows that there exists some (unique) $m$ with $P(x, n, m)$. Combining both facts together with the assumption $H_\varphi(X, Y)$, we conclude that $(x, (n, i)) \in Y$ holds if and only if $seq(m)_i$ is equal to~$1$ for any $i \in \N$. Consequently, we have $Y_{x, n} \mapsto m$ and, therefore, with $g(x, n) = m$ the contradiction $P(x, n, g(x, n))$.
\end{proof}

With this proposition, the proof of Theorem \ref{thm:WETR} is immediate. Moreover, we have the following corollary:

\begin{corollary}[$\RCAn$]
	The following are equivalent:
	\begin{enumerate}[label=\alph*)]
		\item Weak Effective Transfinite Recursion \emph{for all well-orders}
		\item Weak Effective Transfinite Recursion \emph{for well-founded relations}
		\item Weak K\H{o}nig's Lemma
	\end{enumerate}
\end{corollary}

\begin{proof}
	From Proposition \ref{prop:WETR_proves_WKL}, we have that $\WETR$ implies the disjunction of $\WKL$ and $\SETR$. Following Corollary \ref{cor:induction_to_SETR_X} and the fact that $\ACAn$ is a stronger system than $\WKLn$, we conclude that $\WETR$ implies $\WKL$. Finally, Proposition \ref{prop:WKL_proves_WETR'} and the trivial fact that $\WETR'$ has weaker requirements on $X$ than $\WETR$ provide the other two directions.
\end{proof}

We conclude this section with a short discussion on the strength of $\WETR_\omega$: From the previous results, it is clear that this principle is equivalent to the disjunction of $\WKL$ and $\Pi^0_2$-induction along $\N$. At first, it might be surprising that such a natural principle is equivalent to such a seemingly artificial formula. However, the disjunction of $\WKL$ and $\Pi^0_2$-induction along $\N$ has already appeared in the literature: In \cite{Friedman93}, Friedman, Simpson, and Yu discussed this principle in Section~4. Moreover, their Theorem~4.5 states the equivalence of this disjunction with the (very natural) fact that $f^k$ is continuous for all $k \in \N$ if $f: X \to X$ is continuous where $X$ is a metric space. Moreover, in \cite[Theorem 2.24]{Belanger15}, Belanger proves several equivalence results of our disjunction and the implications of certain amalgamation properties. From Theorem 2.16 of the same paper, we see that the disjunction of $\WKL$ and the $\Sigma^0_2$-bounding principle is unprovable in $\RCAn$. With the fact that $\Sigma^0_2$-induction implies this bounding principle (cf.~\cite{Parsons} also \cite[Proposition~4]{KP78}), this implies the unprovability of our disjunction and, therefore, of $\WETR_\omega$ in $\RCAn$.

\subsection*{Acknowledgements} I would like to thank everyone that I and Anton Freund, my PhD advisor, corresponded with during my work on this paper: Lev Beklemishev, Jeffry Hirst, Paul Shafer, and Keita Yokoyama. This, of course, extends to the reviewers, whose suggestions were very helpful. Finally, I would also like to thank Anton Freund for his input and the many discussions that improved both my arguments and the way they were written down. This work was funded by the Deutsche Forschungsgemeinschaft (DFG, German Research Foundation) -- Project number 460597863.

\bibliographystyle{elsarticle-num}
\bibliography{etr}
	
\end{document}